\documentclass[reqno,twoside,11pt,english]{amsart}
\usepackage{amsmath,amsfonts,amssymb,amsthm,epsfig}
\newcommand{\RNum}[1]{\uppercase\expandafter{\romannumeral #1\relax}}
\usepackage{comment}

\usepackage{color}
\usepackage[final,colorlinks,linkcolor=blue,anchorcolor=red,citecolor=blue]{hyperref}
\usepackage{graphicx}
\usepackage{titletoc,epsf}
\usepackage{amsmath,amsfonts,latexsym,amsthm,amsxtra,amssymb,bbm}
\allowdisplaybreaks
\usepackage{graphicx,float}
\usepackage{tikz}
\usetikzlibrary{intersections,patterns}
\usetikzlibrary{calc,spy}
\usepackage{ifthen}
\usepackage{float}

\allowdisplaybreaks

\newtheorem{thm}{Theorem}[section]
\newtheorem{lem}[thm]{Lemma}
\newtheorem{op}[thm]{Open Problem}
\newtheorem{cor}[thm]{Corollary}
\newtheorem{prop}[thm]{Proposition}
\newtheorem{assertion}{Assertion}[section]

\newtheorem{step}{Step}[section]

\newtheorem{cl}{Claim}[section]

\newtheorem{ca}{Case}
\newtheorem{sca}[section]{Subcase}
\newtheorem{scl}[section]{Subclaim}
\newtheorem{conj}[equation]{Conjecture}

\theoremstyle{definition}
\newtheorem{defn}[thm]{Definition}

\newtheorem{ques}[equation]{Question}
\newtheorem{rem}[thm]{Remark}
\newtheorem{exam}[thm]{Example}

\newcounter {own}
\def\theown {\thesection       .\arabic{own}}
\numberwithin{equation}{section}

\newenvironment{pf}[1][]{%
	\vskip 3mm
	\noindent
	\ifthenelse{\equal{#1}{}}%
	{{\slshape Proof. }}%
	{{\slshape #1.} }%
}%
{\qed\bigskip}

\newtheorem{Thm}{Theorem}

\newtheorem{Ques}[Thm]{Question}


\newcommand{\IR}{{\mathbb R}}

\newcommand{\diam}{{\operatorname{diam}}}

\newcommand{\bas}{\begin{assertion}}
	\newcommand{\eas}{\end{assertion}}
\newcommand{\ben}{\begin{enumerate}}
	\newcommand{\een}{\end{enumerate}}
\newcommand{\bst}{\begin{step}}
	\newcommand{\est}{\end{step}}

\def\be{\begin{equation}}
	\def\ee{\end{equation}}

\newcommand{\bee}{\begin{enumerate}}
	\newcommand{\eee}{\end{enumerate}}

\newcommand{\blem}{\begin{lem}}
	\newcommand{\elem}{\end{lem}}
\newcommand{\bthm}{\begin{thm}}
	\newcommand{\ethm}{\end{thm}}
\newcommand{\bcor}{\begin{cor}}
	\newcommand{\ecor}{\end{cor}}
\newcommand{\beg}{\begin{exam}}
	\newcommand{\eeg}{\end{exam}}
\newcommand{\begs}{\begin{examples}}
	\newcommand{\eegs}{\end{examples}}
\newcommand{\bdefe}{\begin{defn}}
	\newcommand{\edefe}{\end{defn}}
\newcommand{\bprob}{\begin{prob}}
	\newcommand{\eprob}{\end{prob}}
\newcommand{\bques}{\begin{ques}}
	\newcommand{\eques}{\end{ques}}
\newcommand{\bei}{\begin{itemize}}
	\newcommand{\eei}{\end{itemize}}
\newcommand{\bcon}{\begin{conj}}
	\newcommand{\econ}{\end{conj}}
\newcommand{\bop}{\begin{op}}
	\newcommand{\eop}{\end{op}}

\newcommand{\bstep}{\begin{step}}
	\newcommand{\estep}{\end{step}}

\newcommand{\bca}{\begin{ca}}
	\newcommand{\eca}{\end{ca}}
\newcommand{\bsca}{\begin{sca}}
	\newcommand{\esca}{\end{sca}}

\newcommand{\bcl}{\begin{cl}}
	\newcommand{\ecl}{\end{cl}}

\newcommand{\bscl}{\begin{scl}}
	\newcommand{\escl}{\end{scl}}

\newcommand{\bcons}{\begin{conjs}}
	\newcommand{\econs}{\end{conjs}}
\newcommand{\bprop}{\begin{prop}}
	\newcommand{\eprop}{\end{prop}}
\newcommand{\br}{\begin{rem}}
	\newcommand{\er}{\end{rem}}
\newcommand{\brs}{\begin{rems}}
	\newcommand{\ers}{\end{rems}}
\newcommand{\bo}{\begin{obser}}
	\newcommand{\eo}{\end{obser}}
\newcommand{\bos}{\begin{obsers}}
	\newcommand{\eos}{\end{obsers}}
\newcommand{\bpf}{\begin{pf}}
	\newcommand{\epf}{\end{pf}}
\newcommand{\ba}{\begin{array}}
	\newcommand{\ea}{\end{array}}
\newcommand{\beq}{\begin{eqnarray}}
	\newcommand{\beqq}{\begin{eqnarray*}}
		\newcommand{\eeq}{\end{eqnarray}}
	\newcommand{\eeqq}{\end{eqnarray*}}

\newcommand{\ds}{\displaystyle}

\newcounter{minutes}
\divide\time by 60
\newcounter{hours}
\multiply\time by 60 \addtocounter{minutes}{-\time}

\begin{document}

	\title[Characterizations of quasihyperbolic John domains and uniform domains]{Characterizations of quasihyperbolic John domains and uniform domains in metric spaces}
	
	\author[S.-J. Gao, C.-Y. Guo, M. Huang and X. Wang]{Shu-Jing Gao, Chang-Yu Guo, Manzi Huang and Xiantao Wang}
	
	\address[S.-J. Gao]{School of Mathematics, Shandong University, Jinan and Research Center for Mathematics and Interdisciplinary Sciences, Shandong University, Qingdao, P. R. China}
	\email{gsjing9@163.com}
	
	\address[C.-Y. Guo]{Research Center for Mathematics and Interdisciplinary Sciences, Shandong University, 266237, Qingdao, P. R. China, and Department of Physics and Mathematics, University of Eastern Finland, 80101, Joensuu, Finland}
	\email{changyu.guo@sdu.edu.cn}
	
	\address[M. Huang]{MOE-LCSM, School of Mathematics and Statistics, Hunan Normal University, Changsha, Hunan 410081, P. R. China} \email{mzhuang@hunnu.edu.cn}
	
	\address[X. Wang]{MOE-LCSM, School of Mathematics and Statistics, Hunan Normal University, Changsha, Hunan 410081, P. R. China} \email{xtwang@hunnu.edu.cn}
	
	\date{}

\begin{abstract}
	In a recent work of Zhou and Ponnusamy [Ann. Sc. Norm. Super. Pisa Ci. Sci. 2025], the authors studied the following natural question:  find sufficient and necessary conditions for a domain $\Omega$ in a metric space $X$ to be quasihyperbolic John. It was proved that Gromov hyperbolic John domains are quasihyperbolic John, quantitatively. As an application, they obtained a characterization of uniform domains in Ahlfors regular spaces. 
	
	In a recent work, using a deep improved characterization of Gromov hyperbolicity, Guo, Huang and Wang [arXiv 2025] proved the quantitative equivalence bteween inner uniformity and the quasihyperbolic John condition in metric doubling spaces. However, the proof does not yield a similar characterization for uniform domains. In this article, we find a new elementary approach to successfully extend the above characterization to uniform domains: a domain $\Omega$ in a doubling length space $X$ is uniform if and only if it is linearly locally connected (LLC) and satisfies the ball separation condition, if and only if it is LLC-1 and quasihyperbolic John, quantitatively. This substantially improved the corresponding results of Zhou and Ponnusamy. Our new approach also allows us to give an alternative proof of the inner uniformity result of Guo-Huang-Wang without using the improved characterization on Gromov hyperbolicity. 
\end{abstract}


\thanks{S.-J. Gao and C.-Y. Guo are supported by the Young Scientist Program of the Ministry of Science and Technology of China (No.~2021YFA1002200), the NSF of China (12311530037), the Taishan Scholar Project and the Jiangsu Provincial Scientific Research Center of Applied Mathematics under Grant No.~BK20233002. M. Huang and X. Wang are partly supported by NSF of China (No.12371071, No.12271189 and No.12571081). }

\maketitle

\section{Instruction}\label{sec-1}


In this paper, we mainly concern geometric characterizations of uniform/inner uniform domains in metric spaces. Recall that 

\begin{defn}\label{def:uniform domain}
	A domain $\Omega$ in a metric space $X=(X,d)$ is called  {\it $c$-uniform}, $c\geq 1$, if each pair of points $z_{1},z_{2}$ in $\Omega$ can be joined by a rectifiable curve $\gamma$ in $\Omega$ satisfying
	\begin{enumerate}
		\item\label{con1} $\ds\min_{j=1,2}\{\ell_d (\gamma [z_j, z])\}\leq c\, d_\Omega(z)$ for all $z\in \gamma$, and
		\item\label{con2} $\ell_d(\gamma)\leq c\,d(z_{1},z_{2})$,
	\end{enumerate}
	where $\ell_d(\gamma)$ denotes the arc-length of $\gamma$ with respect to the metric $d$,
	$\gamma[z_{j},z]$ the subcurve of $\gamma$ between $z_{j}$ and $z$, and $d_\Omega(z):=d(z,\partial \Omega)$ denotes the distance from $z$ to the boundary of $\Omega$. In a $c$-uniform domain $\Omega$, any curve $\gamma\subset \Omega$, which satisfies conditions (1) and (2) above, is called a {\it $c$-uniform curve}.
\end{defn}

If the condition \eqref{con2} in Definition \ref{def:uniform domain} is replaced by the weaker inequality
\begin{equation}\label{eq:def for inner quasiconvexity}
	\ell_d(\gamma)\leq c\,\sigma_\Omega( z_{1}, z_{2}),
\end{equation}
where $\sigma_\Omega$ is the {\it inner distance} in $\Omega$ defined by
$$\sigma_\Omega(z_1,z_2)=\inf \{\ell(\alpha):\; \alpha\subset \Omega\;
\mbox{is a rectifiable curve joining}\; z_1\; \mbox{and}\; z_2 \},$$
then $\Omega$ is said to be $c$-{\it inner uniform} and the corresponding curve $\gamma$ is called a $c$-\emph{inner uniform curve}. When the context is clear, we often drop the subscript $\Omega$ from $\sigma_{\Omega}$ and simply write $\sigma$.

If $\Omega$ only satisfies the condition \eqref{con1} in Definition \ref{def:uniform domain},  then it is said to be a {\it $c$-John domain}, and the corresponding curve $\gamma$ is called a \emph{$c$-John curve} or a {\it double $c$-cone curve}. It is immediate from the definition that 
\[
\{c\text{-uniform domains}\} \subset \{c\text{-inner uniform domains}\}\subset  \{c\text{-John domains}\}.
\]
It is not difficult to construct examples to see that the above inclusions are strict.

The class of John domains in $\mathbb{R}^n$ was initially introduced by F. John in his study of elasticity \cite{Jo} and the name was coined by Martio and Sarvas in \cite{MS}, where they also introduced the class of uniform domains. These classes of domains are central in modern geometric function theory and have wide connections with many other mathematical subjects related to analysis and geometry; see for instance \cite{BK1,GM,GePa,GGGKN-2024,Guo-2015,Hajlasz-Koskela-2000,H,Hei-Book-2001,HeKo,HLPW,Jo81,Martin-1985}.

In a seminal work \cite{BHK}, Bonk-Heinonen-Koskela  established a rich uniformization theory for \emph{Gromov hyperbolic spaces}, using the class of uniform domains as model domains:
\[
\{\text{some class of Gromov hyperbolic spaces}\}=\mathcal{F}\left(\{\text{bounded locally compact uniform spaces}\}\right),
\] 	
where $\mathcal{F}$ is a class of good homeomorphisms (i.e. quasiisometries). To introduce the class of Gromov hyperbolic domains, we first recall the notion of Gromov hyperbolic spaces introduced by M. Gromov in his celebrated work \cite{Gr1}. 
\begin{defn}
	A geodesic metric space $X=(X,d)$ is called \emph{$\delta$-Gromov hyperbolic}, $\delta>0$, if each side of a geodesic triangle in $X$ lies in the $\delta$-neighborhood of the other two sides.
\end{defn}
The Gromov hyperbolicity is a large-scale property, which generalizes the metric properties of classical hyperbolic geometry and of trees, and it turns out to be very useful in geometric group theory and metric geometry \cite{Bonk-Kleiner-2005,BS,Bridson-H-1999,Gr2}. For more about the Gromov hyperbolicity and its connection to geometric function theory, see for instance \cite{BB2,BB4,BHK,Herron-Koskela-1991,Herron-S-Xie-2008,HRWZ-2025}.

Next, we recall the definition of quasihyperbolic metric, which was initially introduced by Gehring and Palka \cite{GePa} for domains in $\IR^n$, and then has been extensively studied in \cite{GO}.
The {\it quasihyperbolic length} of a rectifiable arc
$\gamma$ 
in a proper domain $\Omega\subsetneq (X,d)$ is defined as
$$
\ell_{k}(\gamma):=\int_{\gamma}\frac{|dz|}{d_\Omega(z)}.
$$
For any $z_1$, $z_2$ in $\Omega$, the {\it quasihyperbolic distance}
$k_\Omega(z_1,z_2)$ between $z_1$ and $z_2$ is set to be
$$k_\Omega(z_1,z_2)=\inf_{\gamma}\{\ell_k(\gamma)\},
$$
where the infimum is taken over all rectifiable arcs $\gamma$
joining $z_1$ and $z_2$ in $\Omega$. An arc $\gamma$ from $z_1$ to $z_2$ is called a {\it quasihyperbolic geodesic} if
$\ell_k(\gamma)=k_\Omega(z_1,z_2)$. Clearly, each subarc of a quasihyperbolic
geodesic is a quasihyperbolic geodesic.

\begin{defn}\label{def:Gromov hyperbolic domains}
	Let $(X,d)$ be a metric space and $k$ be 
	the quasihyperbolic distance  {induced by $d$}. A domain $\Omega\subsetneq (X,d)$ is called $\delta$-Gromov hyperbolic if the metric space $(\Omega,k)$ is $\delta$-Gromov hyperbolic.
\end{defn}

In their second main result \cite[Theorem 1.11 and Proposition 7.12]{BHK}, Bonk, Heinonen and Koskela obtained the following useful characterization of uniformity/inner uniformity. To record their result, recall that 
\begin{defn}[LLC]\label{def:LLC}
	Fix $c\geq 1$. A domain $\Omega$ in a metric space $(X,d)$ is said to be $c$-linearly locally connected, $c$-LLC for short, if it satisfies the following two conditions:
	\begin{description}
		\item[$c$-LLC-1] each pair of points $a,b\in \Omega\cap {B(x,r)}$ can be joined in $\Omega\cap {B(x,cr)}$;
		\item[$c$-LLC-2] each pair of points $a,b\in \Omega\backslash \overline{B(x,r)}$ can be joined in $\Omega\backslash \overline{B(x,r/c)}$.
	\end{description}
	It is called $c$-LLC-$i$ ($i\in \{1,2\}$) if only the $i$-th condition is satisfied. 
\end{defn}

With these definitions at hand, the result of Bonk-Heinonen-Koskela can be stated as follows.  
\begin{Thm}[{\cite[Theorem 1.11 and Proposition 7.12]{BHK}}]\label{thm:A}
	A domain $\Omega \subsetneq \mathbb{R}^n$ is $c$-uniform if and only if it is $\delta$-Gromov hyperbolic and $c_0$-LLC, quantitatively. Similarly, a domain $\Omega \subsetneq \mathbb{R}^n$ is $c$-inner uniform if and only if it is $\delta$-Gromov hyperbolic and $c_0$-LLC-2, quantitatively. 
\end{Thm}

By a seminal result of Balogh-Buckley \cite[Theorem 0.1]{BB4}, the class of Gromov hyperbolic domains in $\mathbb{R}^n$ can be characterized by \emph{the Gehring-Hayman inequality} and \emph{the ball separation condition}. Very recently, Guo-Huang-Wang \cite{GHW} proved that Gromov hyperbolic domains in $Q$-doubling length spaces (and thus particularly in $\mathbb{R}^n$) are indeed completely characterized by the geometric ball separation condition. Consequently, Theorem \ref{thm:A} holds if the $\delta$-Gromov hyperbolicity gets replaced by the ball separation condition. Recall that

\begin{defn}\label{def:ball separation}
	A domain $\Omega$ in a metric space $(X,d)$ is said to satisfy the \emph{$C$-ball separation condition}, $C>0$, if for each quasihyperbolic geodesic $\gamma_{xy}\subset \Omega$, every $z\in \gamma_{xy}$ and every curve $\gamma\subset \Omega$ joining $x$ and $y$, it holds
	\begin{equation}\label{eq:def for ball separation}
		B_{\sigma}\left(z,Cd_{\Omega}(z)\right)\cap \gamma\neq \emptyset.
	\end{equation}
	To emphasize the distance $\sigma$, we also say that $(\Omega,\sigma)$ satisfies the $C$-ball separation condition.
\end{defn}

\begin{defn}\label{def:Q doubling ms}
	A metric space $X=(X,d)$ is called \emph{$Q$-doubling}, if for each ball $B(x,r)$, every $r/2$-separated subset of $B(x,r)$ has at most $Q$ points.
\end{defn}
A simple volume comparison implies that $\IR^n$ (equipped with the standard Euclidean distance) is $2^n$-doubling, according to Definition \ref{def:Q doubling ms}.

A central question asked in \cite{ZP-2024-Pisa} is to \emph{characterize quasihyperbolic John domains in metric spaces}. Recall that 
\begin{defn}\label{def:QH John}
	A domain $\Omega$ in a metric space $(X,d)$ is called quasihyperbolic $c$-John if every 
	quasihyperbolic geodesic $\gamma$ in $\Omega$ is a double $c$-cone curve, i.e. $\gamma$ satisfies Definition \ref{def:uniform domain} (1).
\end{defn}

Addressing the above central question, it is important to study the geometry of quasihyperbolic geodesics, as it will be the natural candidate of uniform/inner uniform curves \cite{BHK,GO,Guo-2015,ZLR}.  In \cite[Theorem 1.3]{ZP-2024-Pisa}, the authors proved that \emph{if $\Omega$ is $\delta$-Gromov hyperbolic and $c_0$-John, then it is quasihyperbolic $c$-John, quantitatively}. In a recent work, as an application of a deep improved geometric characterization of Gromov hyperbolicity, the authors proved the following geometric characterization of quasihyperbolic John domains and inner uniform domains, which in particular improves the second assertion in Theorem \ref{thm:A}. 
\begin{Thm}[{\cite[Theorems 6.1 and 6.3]{GHW}}]\label{thm:main result inner uniform}
Let $(X,d)$ be a $Q$-doubling length space and $\Omega$ be a proper subdomain in $X$ such that $(\Omega,k)$ is geodesic. Then the following statements are quantitatively equivalent:
\begin{enumerate}
	\item $\Omega$ is $c$-inner uniform;
	\item $\Omega$ is $c_0$-LLC-2 and satisfies the $c_1$-ball separation condition;
	\item $\Omega$ is quasihyperbolic $c_3$-John.
\end{enumerate}	
\end{Thm}

Note that the assumption of $(\Omega,k)$ being geodesic is quite natural as it is part of the definition of a quasihyperbolic John domain and also the definition of the ball separation condition. If $X$ is a locally compact length space, then the identity map $id:(\Omega,d)\to (\Omega,\sigma_{\Omega})$ is a homeomorphism, and so, by \cite[Proposition 2.8]{BHK}, $(\Omega,k)$ is geodesic. This property remains valid for nice domains in certain infinite dimensional spaces as well; see for instance the work of Martio and V\"{a}is\"{a}l\"{a} \cite{Martio-Vaisala-2011}.  Theorem \ref{thm:main result inner uniform} substantially improves \cite[Theorem 1.3]{ZP-2024-Pisa} above, as it gives a complete characterization of quasihyperbolic John property in terms of inner uniformity or in terms of two geometric condtions (LLC-2 and the ball separation condition).
It also provides a concrete answer to the above central question of Zhou and Ponnusamy in the setting of metric doubling spaces.

As an application of \cite[Theorem 1.3]{ZP-2024-Pisa}, the authors obtained in \cite[Corollary 1.7]{ZP-2024-Pisa} the following characterization of uniformity in Ahlfors regular metric measure spaces: \emph{Let $\Omega$ be minimally nice and Ahlfors $Q$-regular. Then $\Omega$ is uniform if and only if it is quasiconvex, John and roughly starlike Gromov hyperbolic.} We shall not give precise definitions of minimally nice spaces or Ahlfors $Q$-regular spaces/measures, but refers to \cite{ZP-2024-Pisa} for details.  
Here, we would like to point out that the requirement of an Ahlfors regular measure is not natural as all the concepts involved are purely metric. It is thus natural to ask 
\begin{Ques}\label{Ques:Zhou}
Can we obtain a geometric characterization of uniformity, similar to Theorem \ref{thm:main result inner uniform}, in metric spaces without a measure.
\end{Ques}
Note that one can not directly extend the proof of Theorem \ref{thm:main result inner uniform} to uniform domains, as there is no obvious way to verify the the quasiconvexity condition (2) in Definition \ref{def:uniform domain}. New ideas are necessary to find a satisfied solution of Question \ref{Ques:Zhou}.

In this paper, we shall develop an elementary approach to prove the following characterization of uniformity, akin to Theorem \ref{thm:main result inner uniform}, which addresses Question \ref{Ques:Zhou} above.    
\begin{thm}\label{thm:main result uniform}
	Let $(X,d)$ be a $Q$-doubling length space and $\Omega$ be a proper subdomain in $X$ such that $(\Omega,k)$ is geodesic. Then the following statements are quantitatively equivalent:
	\begin{enumerate}
		\item $\Omega$ is $c$-uniform;
		\item $\Omega$ is $c_0$-LLC and satisfies the $c_1$-ball separation condition;
		\item $\Omega$ is $c_2$-LLC-1 and quasihyperbolic $c_3$-John.
	\end{enumerate}
\end{thm}

By \cite[Remark 3.16]{Vai2004}, there is an LLC/LLC-2 Gromov hyperbolic domain $\Omega$ in an infinite dimensional Banach space, which fails to be uniform/inner uniform. This suggests that Theorems \ref{thm:main result uniform} and \ref{thm:main result inner uniform} would probably fail in infinite dimensions. On the other hand, using the new techniques developed in \cite{GHW-1}, Guo-Huang-Li-Wang \cite[Theorem 1.6]{GHLW} did prove that a domain $\Omega$ in a general Banach space $E$ is $c$-inner uniform if it is $\delta$-Gromov hyperbolic and $c_0$-John, quantitatively.

Comparing with \cite[Corollary 1.7]{ZP-2024-Pisa}, our assumptions in Theorem \ref{thm:main result uniform} are much weaker, in particular, we do not require a measure. On the other hand, our conclusions are indeed stronger: LLC and ball separation condition are sufficient to achieve uniformity, while in \cite{ZP-2024-Pisa}, it requires quasiconvexity, John and roughly starlike Gromov hyperbolicity.

The implication  ``(2) $\Rightarrow$ (3)" in Theorem \ref{thm:main result uniform}  has already been proved in \cite[Lemma 6.2]{GHW} (see Lemma \ref{Theorem-1.8} below). The main step for the proof of Theorem \ref{thm:main result uniform} (or Theorem \ref{thm:main result inner uniform}) is thus the reverse implication ``(3) $\Rightarrow$ (1)", as the implication ``(1) $\Rightarrow$ (2)" is rather direct from definition.  The main step for this implication lies in Lemma \ref{Theorem-1.8-0}, which gives length comparison for quasihyperbolic geodesics with different end points. A version of this result was proved earlier in \cite{GHW-1} for Gromov hyperbolic domains in Banach spaces. One can check the proof there does not work in our case. The novelty here is that we found a smart construction of points on a given quasihyperbolic geodesic so that one can obtain effective estimates on the quasihyperblic distance between consecutive points,  using merely the metric doubling condition.  A similar idea is used in the final proof as well. As a contrast, we would like to point out that all the earlier known proof for such a characterization (as in \cite{BHK}, \cite{ZP-2024-Pisa}) relies on the uniformization theory of Bonk-Heinonen-Koskela and thus it requires an Ahlfor regular measure (as well as some other technical assumptions).  


\medskip
\textbf{Structure.} The structure of this paper is as follows. Section \ref{sec:elmentary} contains some elementary results that are necessary for later proofs. In Section \ref{sec:proofs}, we prove our main theorems. As some elementary results from \cite{GHW} are used in our proofs, we add an appendix to include the proofs of these auxiliary results so that it is self-contained and easier for the readers to understand the complete proof.
\smallskip

\textbf{Notations.} Throughout this paper, $X$ is always assumed to be a metric space. For a given metric $d$ on $X$,  $x\in X$ and $r>0$, we define 
$$B_d(x,r):=\{y\in X: d(x,y)<r\}.$$
For each proper subdomain $\Omega\subset X$, we use $\Lambda_{xy}(\Omega)$ to represent the set of all quasihyperbolic geodesics in $\Omega$ with end points $x$ and $y$, and use $\gamma_{xy}$ to denote some quasihyperbolic geodesic in $\Lambda_{xy}(\Omega)$.
Meanwhile, we use  $\Gamma_{xy}(\Omega)$ to represent the set of all curves in $\Omega$ with end points $x$ and $y$.

\section{Some elementary estimates}\label{sec:elmentary}

For any $z_1$, $z_2$ in $\Omega$, let $\gamma\subset \Omega$ be an arc with end points $z_1$ and $z_2$. Then we have the following elementary estimates (see for instance \cite[Section 2]{Vai3})
\beq\label{(2.1)}
\ell_{k}(\gamma)\geq
\log\Big(1+\frac{\ell(\gamma)}{\min\{d_\Omega(z_1), d_\Omega(z_2)\}}\Big)
\eeq
and
\beq\label{(2.2)}
\begin{aligned}
	k_{\Omega}(z_1, z_2) &\geq  \log\Big(1+\frac{\sigma_\Omega(z_1,z_2)}{\min\{d_\Omega(z_1), d_\Omega(z_2)\}}\Big)
	\\
	&\geq
	\log\Big(1+\frac{|z_1-z_2|}{\min\{d_\Omega(z_1), d_\Omega(z_2)\}}\Big)
	\geq
	\Big|\log \frac{d_\Omega(z_2)}{d_\Omega(z_1)}\Big|.
\end{aligned}
\eeq

The follow two elementary lemmas were proved in \cite{GHW}. For convenience, we also include the proofs in the appendix. 

\begin{lem}[{\cite[Lemma 2.1]{GHW}}]\label{lem-3-1}
	Suppose $\Omega\subset X$ is a domain, $u$, $v\in \Omega$, and $\alpha$ is a rectifiable curve in $\Omega$ joining $u$ and $v$. Fix $c\geq 1$. If for each $w\in\alpha$,  $\ell(\alpha[u,w])\leq c d_\Omega(w)$, then
	$$k(u,w)\leq 2c \log\Big(1+\frac{2\ell(\alpha[u,w])}{d_\Omega(u)}\Big)\leq 4c \log\Big(1+\frac{\ell(\alpha[u,w])}{d_\Omega(u)}\Big).$$
\end{lem}

\begin{lem}[{\cite[Lemma 2.2]{GHW}}]\label{Lemma-2.1}
	Let $(X, d)$ be a length space and $\Omega$ a proper domain of $X$. Suppose that there are a constant $a>1$ and two points $x_1$ and $x_2$ in $\Omega$ such that $d(x_1,x_2)\leq a^{-1}d_{\Omega}(x_1)$. Then $$k_{\Omega}(x_1,x_2)\leq \frac{9a}{10(a-1)}\frac{d(x_1,x_2)}{d_{\Omega}(x_1)}\leq \frac{10}{9}(a-1)^{-1}$$
	and $$\ell([x_1,x_2])\leq \frac{10a}{9(a-1)}e^{\frac{10}{9}(a-1)^{-1}}d(x_1,x_2).$$
\end{lem}

We shall repeatedly use the following elementary estimate in $Q$-doubling metric spaces.
\begin{lem}[{\cite[Lemma 4.1.11]{HKST-2015}}]\label{qs-5}
	Let $X=(X,d)$ be a $Q$-doubling metric space with constant $Q$ and $\Omega\subset X$ a domain.	Fix $R>0$ and $a\geq 1$ and let $r=\frac{R}{a}$. Then for any $x\in \Omega$, the ball $B(x,R)$ contains
	at most $b$ balls with radius $r$ such that they are disjoint from each other, where $b\leq Q^{[\log_2 a]}$. Here and hereafter, $[\cdot]$ means the greatest integer part.
\end{lem}

\section{Proofs of the main results}\label{sec:proofs}


In this section, we shall prove Theorem \ref{thm:main result uniform} and Theorem \ref{thm:main result inner uniform}.  It is clear that for Theorem \ref{thm:main result uniform}, it suffices to show  that
\begin{itemize}
	\item[i)] For $x,y \in \Omega$, $\gamma_{xy}\in \Lambda_{xy}(\Omega)$, it holds
	\be\label{H25-0520-1}\ell(\gamma_{xy})\leq e^{19\kappa_0\kappa_1\kappa_2}d(x,y).
	\ee
	
	\item[ii)] For each $z\in\gamma_{xy}$, it holds \be\label{H25-0520-2}\min\left\{\ell(\gamma_{xy}[x,z]),\ell(\gamma_{xy}[y,z])\right\}\leq \kappa_0 d_{\Omega}(z),\ee
	where $\kappa_0=(28ec_0c_1)^2([Q^{2\log_2 70ec_0c_1}]+1)$, $\kappa_1=\kappa_0([Q^{\log_2 15\kappa_0c_0}]+1)$ and $\kappa_2=\kappa_1^{4\kappa_1^2+1}$.
\end{itemize}
Then it follows from i) and ii) that for each $x,y\in \Omega$, each $\gamma_{xy}\in \Lambda_{xy}(\Omega)$ is a $C$-uniform curve and in particular $\Omega$ is $C$-uniform with $C=e^{19\kappa_0\kappa_1\kappa_2}$.

The quasihyperbolic John condition ii) as above has been proved in \cite[Lemma 6.2]{GHW}. For the convenience of readers, we shall repeat its proof in the appendix. 
\blem[{\cite[Lemma 6.2]{GHW}}]\label{Theorem-1.8}
If $\Omega$ is $c_0$-LLC-2 and satisfies the $c_1$-ball separation condition, then it is quasihyperbolic $\kappa_0$-John. More precisely, let $x_1, x_2\in \Omega$ and $\gamma=\gamma_{x_1x_2}\in\Lambda_{x_1x_2}(\Omega)$. Then for each $x\in \gamma$, we have
$$\min\left\{\ell(\gamma[x_1,x]),\ell(\gamma[x_2,x])\right\}\leq \kappa_0 d_{\Omega}(x).$$
\elem

\medskip


\subsection{A technical lemma}

For the proof of statement i), we need the following technical lemma. A version of this lemma was proved in \cite[Lemma 2.27]{GHW-1} for Gromov hyperbolic domains (in Banach spaces).
\blem\label{Theorem-1.8-0} 
Suppose that $\Omega$ is $c_0$-LLC-2 and satisfies the $c_1$-ball separation condition. Let $x$, $y$ and $z$ be three points in $\Omega$, and $\min\{k_{\Omega}(x,y), k_{\Omega}(x,z)\}\geq 1$. Fix $\gamma_{xy}\in \Lambda_{xy}(\Omega)$ and $\gamma_{xz}\in \Lambda_{xz}(\Omega)$. If $k_{\Omega}(y,z)\leq 1$, then $$\max\{\ell(\gamma_{xy}),\ell(\gamma_{xz})\}< e^{9\kappa_0\kappa_2} \min\{\ell(\gamma_{xy}),\ell(\gamma_{xz})\}.$$
\elem
\bpf
It follows from Lemma \ref{Lemma-2.1} and the assumption ``$\min\{k_{\Omega}(x,y), k_{\Omega}(x,z)\}\geq 1$" that
\be\label{H25-0525-1}
d(x,y)\geq \frac{9}{19}\max\{d_{\Omega}(x),d_{\Omega}(y)\}\;\mbox{ and }\; d(x,z)\geq \frac{9}{19}\max\{d_{\Omega}(x),d_{\Omega}(z)\}.
\ee

Without loss of generality, we may assume that $$\min\{\ell(\gamma_{xy}),\ell(\gamma_{xz})\}=\ell(\gamma_{xy}).$$
Note that by \eqref{(2.2)}, we have $$\max\left\{\log\Big(1+\frac{\sigma(y,z)}{\min\{d_{\Omega}(y),d_{\Omega}(z)\}}\Big),\ \log \frac{\max\{d_{\Omega}(y),d_{\Omega}(z)\}}{\min\{d_{\Omega}(y),d_{\Omega}(z)\}}\right\}\leq k_{\Omega}(y,z)\leq 1.$$
It follows that
\be\label{H25-0514-1}\sigma(y,z)\leq (e-1)\min\{d_{\Omega}(y),d_{\Omega}(z)\}\ee and
\be\label{H25-0514-2} \max\{d_{\Omega}(y),d_{\Omega}(z)\}\leq e \min\{d_{\Omega}(y),d_{\Omega}(z)\}.\ee

Take $y_0\in\gamma_{xy}$ and $z_0\in\gamma_{xz}$ such that
\be\label{H25-0514-3}\ell(\gamma_{xy}[x,y_0])=\frac{1}{2}\ell(\gamma_{xy})\;\mbox{ and }\;\ell(\gamma_{xz}[z,z_0])=\frac{1}{2}\ell(\gamma_{xz}).
\ee
We consider the following two cases.

\medskip

\textbf{Case I:} $\ell(\gamma_{xy})\leq e^{\kappa_2}\min\{d_{\Omega}(y),d_{\Omega}(x)\}$.
\medskip

In this case, we have
\begin{equation}\label{eq:lemma 6.2}
	\log \Big(1+\frac{\ell(\gamma_{xz})}{\min\{d_{\Omega}(z),d_{\Omega}(x)\}}\Big)\stackrel{\eqref{(2.1)}}{\leq} k_{\Omega}(x,z)\leq k_{\Omega}(x,y)+k(y,z)\leq k_{\Omega}(x,y_0)+k_{\Omega}(y,y_0)+1.
\end{equation}
Note that by Lemma \ref{Theorem-1.8} and the choice of $y_0$, it holds
$$
\ell(\gamma_{xy}[x,y_0])\leq \kappa_0 d_{\Omega}(y_0) \quad \text{and}\quad \ell(\gamma_{xy}[y,y_0])\leq \kappa_0 d_{\Omega}(y_0).
$$
Thus it follows from Lemma \ref{lem-3-1} that
\[
k_{\Omega}(x,y_0)\leq 4\kappa_0\log \Big(1+\frac{\ell(\gamma_{xy}[x,y_0])}{d_{\Omega}(x)}\Big)
\]
and 
\[
k_{\Omega}(y,y_0)\leq 4\kappa_0\log \Big(1+\frac{\ell(\gamma_{xy}[y,y_0])}{d_{\Omega}(y)}\Big).
\]
Substituting these estimates in \eqref{eq:lemma 6.2}, we obtain that
\beqq
\begin{aligned}
	\log \Big(1+\frac{\ell(\gamma_{xz})}{\min\{d_{\Omega}(z),d_{\Omega}(x)\}}\Big) &\leq 4\kappa_0\left(\log \Big(1+\frac{\ell(\gamma_{xy}[x,y_0])}{d_{\Omega}(x)}\Big)+\log \Big(1+\frac{\ell(\gamma_{xy}[y,y_0])}{d_{\Omega}(y)}\Big)\right)+1
	\\
	&< 8\kappa_0\kappa_2 \stackrel{\eqref{H25-0525-1}}{<} \frac{20\kappa_0\kappa_2\ell(\gamma_{xy})}{\min\{d_{\Omega}(x),d_{\Omega}(y)\}} \stackrel{\eqref{H25-0514-2}}{\le}\frac{20e\kappa_0\kappa_2\ell(\gamma_{xy})}{\min\{d_{\Omega}(x),d_{\Omega}(z)\}}.
\end{aligned}
\eeqq
Since $\log(1+t)\geq a_1t$ for $t\leq a_1^{-1}-1$ and $a_1>0$, it follows from the above inequality that
$$  { \frac{e^{-8\kappa_0\kappa_2}\ell(\gamma_{xz})}{\min\{d_{\Omega}(x),d_{\Omega}(z)\}}\leq \log \Big(1+\frac{\ell(\gamma_{xz})}{\min\{d_{\Omega}(z),d_{\Omega}(x)\}}\Big)}\leq 20e\kappa_0\kappa_2\frac{\ell(\gamma_{xy})}{\min\{d_{\Omega}(x),d_{\Omega}(z)\}}$$
and so
\be\label{H25-0515-2}
\ell(\gamma_{xz})\leq 20\kappa_0\kappa_2e^{1+8\kappa_0\kappa_2}\ell(\gamma_{xy})\leq e^{9\kappa_0\kappa_2}\min\{\ell(\gamma_{xy}),\ell(\gamma_{xz})\}.
\ee
This completes the proof  for this case.

\medskip

\textbf{Case II:} $\ell(\gamma_{xy})> e^{\kappa_2}\min\{d_{\Omega}(x),d_{\Omega}(y)\}$.
\medskip

In this case, select $m_1\geq[\kappa_2]-1$ such that
$$e^{m_1+1}d_{\Omega}(y)<{\ell(\gamma_{xy})}\leq e^{m_1+2}d_{\Omega}(y).$$
No loss of generality, by (\ref{H25-0514-2}), we may assume that
$$\min\{d_{\Omega}(x),d_{\Omega}(y)\}=d_{\Omega}(y),$$
since the proof for the other case is similar.

Set $y_1=y$ and $z_1=z$. { Choose $y_{m_1+1} \in\gamma_{xy_1}[y_1,y_0]$, $y_{1,m_1+1} \in\gamma_{xy_1}[x,y_0]$, $z_{m_1+1} \in \gamma_{xz_1}[z_1,z_0]$ and $z_{1,m_1+1} \in \gamma_{xz_1}[x,z_0]$ such that
	\beq\label{5-14-1}
	\begin{aligned} 
		\ell (\gamma_{xy_1}[y_1,y_{m_1+1}]) &=\ell (\gamma_{xz_1}[z_1,z_{m_1+1}])=\ell (\gamma_{xy_1}[x,y_{1,m_1+1}])=\ell (\gamma_{xz_1}[x,z_{1,m_1+1}])\\
		&=e^{m_1}d_{\Omega}(y_1).
	\end{aligned}
	\eeq
}

Next, we claim that
\bcl\label{H25-0520}$\max \{k_{\Omega}(y_{m_1+1},z_{m_1+1}), k_{\Omega}(y_{1,m_1+1},z_{1,m_1+1})\} \leq 2\kappa_2.$\ecl

Suppose on the contrary that  $$\max \{k_{\Omega}(y_{m_1+1},z_{m_1+1}), k_{\Omega}(y_{1,m_1+1},z_{1,m_1+1})\} > 2\kappa_2.$$
Without loss of generality, we may assume that
\beq\label{5-14-3}
k_{\Omega}(y_{m_1+1},z_{m_1+1}) >2\kappa_2,
\eeq
as the proof for the other case is similar.

For each $i \in \{2,\cdots,m_1+1\}$, let $y_i \in \gamma_{xy_1}[y_1 ,y_0 ]$ and $z_i \in \gamma_{xz_1}[z_1 ,z_0 ]$  be such that
\beq\label{5-14-4}
\ell (\gamma_{xy_1}[y_1 ,y_i ]) =\ell (\gamma_{xz_1}[z_1 ,z_i ])=e^{i-1} d_{\Omega}(y_{1}).
\eeq
Then, for each  $i \in \{1,\cdots,m_1 \}$, we have
\[
\ell(\gamma_{xy_1}[y_{i+1},y_i])=\ell(\gamma_{xy_1}[y_1,y_{i+1}])-\ell(\gamma_{xy_1}[y_1,y_{i}])\stackrel{\eqref{5-14-4}}{=}(e-1)e^{i-1}d_{\Omega}(y_1).
\]
By Lemma \ref{Theorem-1.8} and our choice of $y_{i+1}$ in \eqref{5-14-4}, it holds
\be\label{H25-0527-1}
\ell(\gamma_{xy_1}[y_1,y_{i+1}])=e^{i} d_{\Omega}(y_{1})\leq \kappa_0 d_{\Omega}(y_{i+1}).
\ee
Thus, combining the above two estimates gives
\[
\ell(\gamma_{xy_1}[y_{i+1}, y_i])\leq (e-1)\kappa_0 d_{\Omega}(y_i)
\]
and
\[
\ell(\gamma_{xy_1}[y_{i+1}, y_i])\leq\frac{e-1}{e}\kappa_0 d_{\Omega}(y_{i+1}).
\]
For each  $i \in \{1,\cdots,m_1 \}$, it follows from the above inequality and Lemma \ref{lem-3-1} that
\beq\label{5-14-7}
k_{\Omega}(y_{i},y_{i+1})\leq 4\kappa_0  \log\Big(1+\frac{\ell (\gamma_{xy_1}[y_i,y_{i+1}])}{d_{\Omega}(y_{i})}\Big)\leq  4\kappa_0 \log(1+e{\kappa_0}).
\eeq
Similarly, one can prove that
\be\label{H-25-0515-7}
k_{\Omega}(z_{i},z_{i+1})\leq 4\kappa_0 \log(1+e{\kappa_0}).
\ee

Based on \eqref{5-14-3}, we take $t\in\{2,\cdots, {m_1+1}\}$ to be the smallest integer such that
$$k_{\Omega}(y_t,z_t)\geq \kappa_2.$$
Then, it follows from \eqref{5-14-7} and \eqref{H-25-0515-7} that
$$
\kappa_2\leq k_{\Omega}(y_t,z_t)\leq \sum_{i=1}^{t-1}\left(k_{\Omega}(y_i,y_{i+1})+k_{\Omega}(z_i,z_{i+1})\right)+k_{\Omega}(y,z)\leq8(t-1)\kappa_0\log(1+e{\kappa_0})+1.
$$
Thus, we obtain
\be\label{H25-0519-1}
t> \frac{\kappa_2}{8\kappa_0\log(1+e{\kappa_0})}.
\ee

Now, let $s_1=t$, and based on \eqref{H25-0519-1}, for each $i\in\{2,\cdots,[\kappa_1^2]\}$, set $s_i=s_{i-1}-[\kappa_1^{2i}]$.
For each $i\in\{2,\cdots, [\kappa_1^2]\}$, by our choices of $\kappa_0$, $\kappa_1$ and $\kappa_2$,  we have
$$\ell(\gamma_{xy_1}[y_1,y_{s_i}])\geq d_{\Omega}(y_{s_i})-d_{\Omega}(y_1)\stackrel{\eqref{H25-0527-1}}{\geq} (1-\kappa_0e^{1-s_i})d_{\Omega}(y_{s_i})>\frac{1}{2}d_{\Omega}(y_{s_i}).$$
Then
\be\label{H25-0519-0}
{\frac{1}{2}e^{[\kappa_1^{2i}]}d_{\Omega}(y_{s_i})\leq e^{[\kappa_1^{2i}]}\ell(\gamma_{xy_1}[y_1,y_{s_i}])}\stackrel{\eqref{5-14-4}}{=}\ell(\gamma_{xy_1}[y_1,y_{s_{i-1}}])\stackrel{\text{Lemma }\ref{Theorem-1.8}}{\leq} \kappa_0d_{\Omega}(y_{s_{i-1}})\ee
and 
$$d_{\Omega}(y_{s_{i-1}})\leq d_{\Omega}(y_1)+\ell(\gamma_{xy_1}[y_1,y_{s_{i-1}}])\stackrel{\eqref{5-14-4}}{<}2e^{[\kappa_1^{2i}]}\ell(\gamma_{xy_1}[y_1,y_{s_i}])\stackrel{\text{Lemma }\ref{Theorem-1.8}}{\leq} 2\kappa_0e^{[\kappa_1^{2i}]}d_{\Omega}(y_{s_i}).$$
It follows from these two estimates that
\be\label{H25-0519-2}
\frac{1}{2\kappa_0}e^{[\kappa_1^{2i}]}d_{\Omega}(y_{s_i})\leq d_{\Omega}(y_{s_{i-1}})\leq 2\kappa_0e^{[\kappa_1^{2i}]}d_{\Omega}(y_{s_i}).\ee
Similarly, one can prove
\be\label{H25-0519-3}
\frac{1}{2}d_{\Omega}(y_{s_i})\leq\ell(\gamma_{xy_1}[y_1,y_{s_i}])=\ell(\gamma_{xz_1}[z_1,z_{s_i}])\leq \kappa_0d_{\Omega}(z_{s_i}).
\ee

For each $i\in\{2,\cdots, [\kappa_1^2]\}$ and $\gamma_{y_{s_i}z_{s_i}}\in \Lambda_{y_{s_i}z_{s_i}}$, let $x_{0,i}$ bisect the length of $\gamma_{y_{s_i}z_{s_i}}$. Then we have
\[
\begin{aligned}
	k_{\Omega}(y_{s_i}, z_{s_i})&\geq k_{\Omega}(y_t,z_t)-k_{\Omega}(y_t,y_{s_i})-k_{\Omega}(z_t,z_{s_i})\\
	&\geq k_{\Omega}(y_t,z_t)-\sum_{l=s_i}^{t-1}\left(k_{\Omega}(y_{l+1},y_l)+k_{\Omega}(z_{l+1},z_l)\right)\\
	&\stackrel{\eqref{5-14-7}+\eqref{H-25-0515-7}}{\geq} \kappa_2-8\kappa_0{(t-s_i-1)}\log(1+e\kappa_0)
	\\ &> \kappa_2-\kappa_1^{2[\kappa_1^2]+2},
\end{aligned}
\]
where in the last inequality, we used the fact that $t-s_i=\sum\limits_{l=2}^i [\kappa_1^{2l}]<2\kappa_1^{2i}$.

On the other hand, it follows from Lemmas  \ref{lem-3-1} and \ref{Theorem-1.8}, and \eqref{H25-0519-3} that
\beqq
\begin{aligned}
	k_{\Omega}(y_{s_i}, z_{s_i})&= k_{\Omega}(y_{s_i}, x_{0,i})+k_{\Omega}(x_{0,i}, z_{s_i})\\
	&\leq 4\kappa_0\left(\log \left(1+\frac{\ell(\gamma_{y_{s_i}z_{s_i}})}{2d_{\Omega}(y_{s_i})}\right)+\log \left(1+\frac{\ell(\gamma_{y_{s_i}z_{s_i}})}{2d_{\Omega}(z_{s_i})}\right)\right)
	\\ &\leq 4\kappa_0\left(\log (1+\frac{\ell(\gamma_{y_{s_i}z_{s_i}})}{2d_{\Omega}(y_{s_i})})+\log \left(1+\frac{\kappa_0\ell(\gamma_{y_{s_i}z_{s_i}})}{d_{\Omega}(y_{s_i})}\right)\right).
\end{aligned}
\eeqq
Combining the above two estimates with Lemma \ref{Theorem-1.8} gives
\be\label{H25-0519-4}
\ell(\gamma_{y_{s_i}z_{s_i}})>e^{\frac{\kappa_2}{9\kappa_0}} d_{\Omega}(y_{s_i})\geq \frac{1}{\kappa_0}e^{\frac{\kappa_2}{9\kappa_0}}\ell(\gamma_{xy_1}[y_1,y_{s_i}])\geq \frac{1}{\kappa_0}e^{\frac{\kappa_2}{9\kappa_0}-{\color{blue}2}\kappa_1^{2i}}\ell(\gamma_{xy_1}[y_1,y_{s_1}]).\ee


For each $i\in\{2,\cdots, [\kappa_1^2]\}$, based on \eqref{H25-0519-4}, we may select 
$y_{1,i}\in\gamma_{y_{s_i}z_{s_i}}[y_{s_i},x_{0,i}]$ and $z_{1,i}\in\gamma_{y_{s_i}z_{s_i}}[z_{s_i},x_{0,i}]$ so that
\be\label{H25-0519-5}
\ell(\gamma_{y_{s_i}z_{s_i}}[y_{s_i},y_{1,i}])=\ell(\gamma_{y_{s_i}z_{s_i}}[z_{s_i},z_{1,i}])=\ell(\gamma_{xy_1}[y_1,y_{s_1}]).
\ee
Note that
\be\label{H25-0519-7}
k_{\Omega}(z_{s_1}, z_{s_{i-1}})\leq \sum_{l=s_{i-1}}^{s_1-1}k_{\Omega}(z_{l+1}, z_{l})\stackrel{\eqref{H-25-0515-7}}{\leq}8\kappa_0\kappa_1^{2(i-1)}\log(1+e{\kappa_0})
\ee
and
\be\label{H25-0519-7-2}
k_{\Omega}(y_{s_1}, y_{s_{i-1}})\leq \sum_{l=s_{i-1}}^{s_1-1}k_{\Omega}(y_{l+1}, y_{l})\stackrel{\eqref{5-14-7}}{\leq}8\kappa_0\kappa_1^{2(i-1)}\log(1+e{\kappa_0})
\ee
Thus, we obtain from Lemmas  \ref{lem-3-1} and \ref{Theorem-1.8}, \eqref{H25-0527-1} and \eqref{H25-0519-5} that
\[
\begin{aligned}
	k_{\Omega}(y_{s_{i-1}}, y_{1,i-1})&\leq 4\kappa_0\log \Big(1+\frac{\ell(\gamma_{y_{s_{i-1}}z_{s_{i-1}}}[y_{s_{i-1}},y_{1,i-1}])}{d_{\Omega}(y_{s_{i-1}})} \Big)
	\\ &= 4\kappa_0\log \Big(1+\frac{\ell(\gamma_{xy_1}[y_1,y_{s_1}])}{d_{\Omega}(y_{s_{i-1}})}\Big)< 12\kappa_0\kappa_1^{2(i-1)}.
\end{aligned}
\]
On the other hand, 
\[
\begin{aligned}k_{\Omega}(y_{s_i}, y_{1,i})&\stackrel{\eqref{(2.1)}}\geq \log \frac{\ell(\gamma_{y_{s_i}z_{s_i}}[y_{s_i},y_{1,i}])}{d_{\Omega}(y_{s_i})}\stackrel{\eqref{H25-0519-3}+\eqref{H25-0519-5}}\geq
	\log \frac{\ell(\gamma_{xy_1}[y_1,y_{s_1}])}{2\ell(\gamma_{xy_1}[y_1,y_{s_i}])}\\&\stackrel{\eqref{5-14-4}}>\kappa_1^{2i}.
\end{aligned}
\]
The above two inequalities yield that
\be\label{H25-0519-18}k_{\Omega}(y_{s_i}, y_{1,i})>k_{\Omega}(y_{s_{i-1}}, y_{1,i-1})+(\kappa_1-12\kappa_0)\kappa_1^{2i-1}.\ee
A similar discussion as in (\ref{H25-0519-18}) shows that
\be\label{H25-0519-19}k_{\Omega}(z_{s_i}, z_{1,i})>k_{\Omega}(z_{s_{i-1}}, z_{1,i-1})+(\kappa_1-12\kappa_0)\kappa_1^{2i-1}.\ee

Let $B_0=B(y_1, 3\ell(\gamma_{xy_1}[y_1,y_{s_1}]))$. For each $t \in \{2, \dots, [\kappa_1^2]\}$, we take
$$B_{t}^1=B\left(y_{1,t}, \frac{1}{2\kappa_0}\ell(\gamma_{xy_1}[y_1,y_{s_1}])\right)\;\quad\mbox{and}\;\quad B_{1,t}^1=B\left(z_{1,t}, \frac{1}{2\kappa_0}\ell(\gamma_{xy_1}[y_1,y_{s_1}])\right).$$
For each $u \in (\overline{B_t^1}\cup \overline{B_{1,t}^1})$, by (\ref{H25-0519-5}),
$$ d(y_1, u) < 3\ell(\gamma_{xy_1}[y_1,y_{s_1}]),
$$
and so $$ (\overline{B_t^1}\cup \overline{B_{1,t}^1})\subset B_0.$$

Next, we observe that there exist $m=[Q^{\log_2 6\kappa_0}]+1$ balls $B_{1_j}^1$ in $\bigcup\limits_{t=1}^{[\kappa_1^2]}\{B_t^1\}$ such that for each $i\not=j\in\{2,\cdots,m\}$, it holds 
\be\label{H25-10-15-1} 
B_{1_i}^1\cap B_{1_j}^1\not=\emptyset.
\ee
To get \eqref{H25-10-15-1} we come to prove the following claim.
\bcl\label{H25-12-21}For each $t_1<t_2\in\{1,\cdots,[\kappa_1^2]\}$, if $t_2-t_1\geq [\kappa_1^{\frac{1}{2}}]$, there must exist some $s_1\in\{t_1,\cdots,t_2\}$ such that there are at least $[\frac{t_2-t_1}{m^2}]$ balls which satisfy 
\begin{itemize}
	\item\label{H25-12-21-1} for each $j\in\{1,\cdots,[\frac{t_2-t_1}{m^2}]-1\}$, $t_2>q_{j+1}>q_j>s_1$;
	\item\label{H25-12-21-2} and for each $j\in\{1,\cdots,[\frac{t_2-t_1}{m^2}]\}$, $B_{q_j}^1\cap B_{s_1}^1\not=\emptyset$.
\end{itemize}
\ecl
We will prove this claim by a contradiction. Suppose on the contrary that for each $s\in\{t_1,\cdots,t_2\}$, there exist at most $s_0\leq [\frac{t_2-t_1}{m^2}]$ balls $B_{p_j}^1$
such that for each $j\in\{1,\cdots,s_0\}$, 
\be\label{H25-12-21-3}p_j>s\,\mbox{and}\, B_{p_j}^1\cap B_{s}^1\not=\emptyset.\ee
Then (\ref{H25-12-21-3}) implies that there exist at least $[(1-\frac{1}{m^2})(t_2-t_1)]$ balls
$B_{t_{1,j}}^1$ which are disjoint with $B_{t_1}^1$. We take $B_{t_{0,1}}^1=B_{t_1}^1$. Similarly, by (\ref{H25-12-21-3}), there exist at least $[(1-\frac{2}{m^2})(t_2-t_1)]$ balls $B_{t_{2,j}}^1\subset \cup_{j=1}^{[(1-\frac{1}{m})(t_2-t_1)]}B_{t_{1,j}}^1$ which are disjoint with $B_{t_{1,1}}^1$. Repeating the procedure $m$ times we shall get $m$ balls $B_{t_{i,1}}^1$ 
which are disjoint from each other, which contradicts with Lemma \ref{qs-5}. Hence, Claim \ref{H25-12-21} holds.

For each $i\in\{1,\cdots, [\kappa_1^2]\}$, it follows from  Claim \ref{H25-12-21} that there  exist some $1_1\in\{1,\cdots,[\kappa_1^2]\}$ such that there are $n_1\geq [\frac{[\kappa_1^2]}{m^2}]$ balls $B_{q_{1,j}}^1$ which intersect with $B_{1_1}^1$. Similarly, there  exist some $B_{1_2}^1\subset \cup_{j=1}^{n_2}B_{q_{1,j}}^1$ such that there are at least $n_2\geq [\frac{n_1}{m^2}]$ balls $B_{q_{2,j}}^1$ which intersect with $B_{1_2}^1$. Repeating the procedure $m$ times we shall get \eqref{H25-10-15-1}.


It follows from  Lemma \ref{qs-5} and \eqref{H25-10-15-1} that 
there must exist $p<q\in\{2,\cdots, [\kappa_1^2]\}$ such that
\be\label{H25-10-15-2}
B(y_{1,p}, r_2)\cap B(y_{1,q}, r_2)\not=\emptyset\quad\;\mbox{and}\;B(z_{1,p}, r_2)\cap B(z_{1,q}, r_2)\not=\emptyset.
\ee

Let $y_{1,p}^1\in [y_{1,p},y_{1,q}]\cap (B_{p}^1\cap B_q^1)$ and $z_{1,q}^1\in [z_{1,p},z_{1,q}]\cap (B_{1,p}^1\cap B_{1,q}^1)$.
Then, based on Lemma \ref{Theorem-1.8}, \eqref{H25-0519-5}, and \eqref{H25-10-15-2}, we have $$d(y_{1,p}^1,y_{1,p})\leq\frac{1}{2\kappa_0}\ell(\gamma_{xy_1}[y_1,y_{s_1}])=\frac{1}{2\kappa_0}\ell(\gamma_{y_{s_p}z_{s_p}}[y_{s_p},y_{1,p}])\leq\frac{1}{2}d_{\Omega}(y_{1,p}).$$
Similarly, $$d(y_{1,p}^1,y_{1,q})\leq\frac{1}{2}d_{\Omega}(y_{1,q}),\quad\; d(z_{1,p}^1,z_{1,p})\leq\frac{1}{2}d_{\Omega}(z_{1,p})\quad\;\mbox{and}\quad\;d(z_{1,p}^1,z_{1,q})\leq\frac{1}{2}d_{\Omega}(z_{1,q}).$$
Thus we obtain from Lemma \ref{Lemma-2.1} that
$$
\begin{aligned}
	&\quad\max\{k_{\Omega}(y_{1,p},y_{1,q}), k_{\Omega}(z_{1,p},z_{1,q})\}\\
	&\leq \max\{k_{\Omega}(y_{1,p},y_{1,p}^1)+k_{\Omega}(y_{1,q},y_{1,p}^1), k_{\Omega}(z_{1,p},z_{1,p}^1)+k_{\Omega}(z_{1,q},z_{1,p}^1)\}\leq\frac{20}{9},
\end{aligned}
$$
which, together with \eqref{H25-0519-18} and  \eqref{H25-0519-19}, yields
\begin{eqnarray*}k_{\Omega}(y_{s_q}, z_{s_q})\nonumber&=&k_{\Omega}(y_{s_q}, y_{1,q})+k_{\Omega}(y_{1,q}, z_{1,q})+k_{\Omega}(z_{1,q}, z_{s_q}) \\ \nonumber&>&k_{\Omega}(y_{s_p}, y_{1,p})+k_{\Omega}(y_{1,q}, z_{1,q})+k_{\Omega}(z_{1,p}, z_{s_p})+2(\kappa_1-12\kappa_0)\kappa_1^{2q-1}
	\\ \nonumber&\geq&k_{\Omega}(y_{s_p}, y_{1,p})+k_{\Omega}(y_{1,p}, z_{1,p})+k_{\Omega}(z_{1,p}, z_{s_p})+2(\kappa_1-12\kappa_0)\kappa_1^{2q-1}-\frac{40}{9}
	\\ \nonumber&=&k_{\Omega}(y_{s_p},z_{s_p})+2(\kappa_1-12\kappa_0)\kappa_1^{2q-1}-\frac{40}{9},\end{eqnarray*}
and so (\ref{H25-0519-7}) shows that
\begin{eqnarray*}k_{\Omega}(y_{s_q}, z_{s_q})\nonumber&>&k_{\Omega}(y_{s_p},z_{s_p})+k_{\Omega}(y_{s_1},y_{s_p})+k_{\Omega}(z_{s_1},z_{s_p})\\ \nonumber&+&
	2(\kappa_1-12\kappa_0)\kappa_1^{2q-1}-16\kappa_0\kappa_1^{2p}\log(1+e{\kappa_0})-\frac{40}{9}
	\\ \nonumber&>& k_{\Omega}(y_{s_1},z_{s_1})=k_{\Omega}(y_t,z_t),\end{eqnarray*}
which contradicts with the construction
of $y_t$. Hence Claim \ref{H25-0520} holds.
\medskip

By Claim \ref{H25-0520}, we have 
\begin{align} \nonumber\label{5-14-2}
	k_{\Omega}(z_{1,m_1+1},z_{m_1+1}) &\leq k_{\Omega}(y_{m_1+1},y_{1,m_1+1})+k_{\Omega}(y_{1,m_1+1},z_{1,m_1+1})+k_{\Omega}(y_{m_1+1},z_{m_1+1})\\
	&\leq  k_{\Omega}(y_{m_1+1},y_{1,m_1+1})+4\kappa_2.
\end{align}
Meanwhile, we know from Lemma \ref{Theorem-1.8}, Lemma \ref{lem-3-1}, \eqref{H25-0514-3} and \eqref{5-14-4} that
\begin{align} \nonumber
	k_{\Omega}(y_{m_1+1},y_{1,m_1+1}) &= k_{\Omega}(y_{m_1+1},y_{0})+k_{\Omega}(y_{1,m_1+1},y_{0}) \\ \nonumber
	&\leq 4\kappa_0 \log\Big(1+\frac{\ell (\gamma_{xy_1}[y_{m_1+1},y_0])}{d_{\Omega}(y_{m_1+1})}\Big)+4\kappa_0 \log\Big(1+\frac{\ell (\gamma_{xy_1}[y_{1,m_1+1},y_0])}{d_{\Omega}(y_{1,m_1+1})}\Big)\\ \nonumber
	&< 8 \kappa_0 \log(1+e^2\kappa_0),
\end{align}
which, together with \eqref{(2.1)} and \eqref{5-14-2}, shows that
\beqq
\log\Big(1+\frac{\ell (\gamma_{xz_1}[z_{1,m_1+1},z_{m_1+1}])}{d_{\Omega}(z_{m_1+1})}\Big) \leq  k_{\Omega}(z_{1,m_1+1},z_{m_1+1})<5\kappa_2.
\eeqq
Hence, by \eqref{H25-0514-2} and the assumption of Case II, it holds 
\begin{align} \nonumber
	\ell (\gamma_{xz_1}[z_{1,m_1+1},z_{m_1+1}]) &< e^{5\kappa_2}d_{\Omega}(z_{m_1+1})\\ \nonumber
	&<e^{5\kappa_2}(d_{\Omega}(z_{1})+\ell (\gamma_{xy_1}[y_1,y_{m_1+1}]))<e^{5\kappa_2}\ell (\gamma_{xy}).
\end{align}
Finally, we get from (\ref{5-14-1}) that
\begin{align} \nonumber
	\ell (\gamma_{xz_1})&=\ell (\gamma_{xz_1}[x,z_{1,m_1+1}])+\ell (\gamma_{xz_1}[z_{1,m_1+1},z_{m_1+1}])+\ell (\gamma_{xz_1}[y_{1},z_{m_1+1}])\\ \nonumber
	&< (e^{5\kappa_2}+1) \ell (\gamma_{xy}).
\end{align}
The proof of Lemma \ref{Theorem-1.8-0} is thus complete.
\epf

\subsection{Proofs of Theorems \ref{thm:main result uniform} and \ref{thm:main result inner uniform}}

\begin{proof}[Proof of Theorem \ref{thm:main result uniform}]
	(1) $\Rightarrow$ (2): Since $\Omega$ is $C$-uniform, it is  $(2C+1)$-LLC-2 by \cite[Lemma 3.7]{Vai7}. For any $x,y\in B(a,r)\subset\Omega$, let $\gamma$ be a $C$-uniform curve joining $x$ and $y$. For each $z\in \gamma$, we have  
	$$d(a,z)\leq \min\{d(a,x)+d(x,z),d(a,y)+d(y,z)\}\leq r+\frac{1}{2}\ell(\gamma)\leq(2C+1)r,$$ and so $\gamma\subset B(a,(2C+1)r)$. This implies that $\Omega$ is $(2C+1)$-LLC-1, and then $\Omega$ is $(2C+1)$-LLC.  On the other hand, it follows from \cite[Lemma 2.35]{Vai7} (or the proof of \cite[Lemma 2.13]{BHK}) that every quasihyperbolic geodesic $\gamma$ in a $C$-uniform domain $\Omega$ is a $C_1$-uniform curve with $C_1=130C^4e^{32C^4}$, and thus the domain $\Omega$ satisfies the $C_1$-ball separation property with $C_1=130C^4e^{32C^4}$.
	\medskip 
	
	(2) $\Rightarrow$ (3): This follows immediately from Lemma \ref{Theorem-1.8}.  
	
	(3) $\Rightarrow$ (1): For \eqref{H25-0520-1}, we need to show that for any $\gamma_{xy}\in \Lambda_{xy}(\Omega)$, it holds
	$$\ell(\gamma_{xy})\leq e^{19\kappa_0\kappa_1\kappa_2}d(x, y).$$
	
	Suppose on the contrary that there exists some $\gamma_0\in\Lambda_{xy}(\Omega)$  such that
	$$\ell(\gamma_0)> e^{19\kappa_0\kappa_1\kappa_2}d(x, y).$$
	Then
	\be\label{H25-0525-2}
	d(x,y)\ge\frac{1}{2}\max\{d_\Omega(x),d_\Omega(y)\}.\ee
	Otherwise, by Lemma \ref{Lemma-2.1},
	$$\ell(\gamma_0)\le\frac{20}{9}e^{\frac{10}{9}}d(x,y).$$
	This contradiction shows that \eqref{H25-0525-2} holds.
	
	Thus, by \eqref{(2.1)} and \eqref{H25-0525-2}, we have
	$$k_\Omega(x,y)=\ell_k(\gamma_0)\ge \log \Big(1+\frac{\ell(\gamma_0)}{\min\{d_\Omega(x),d_\Omega(y)\}}\Big)>1.$$
	Hence it follows from Lemma \ref{Theorem-1.8-0} that for any $\gamma\in \Lambda_{xy}(\Omega)$,
	\be\label{H25-0525-3}
	\ell(\gamma)\geq  e^{-9\kappa_0\kappa_2}\ell(\gamma_0)> e^{(19\kappa_0-9)\kappa_1\kappa_2}d(x,y).\ee
	
	Since $\Omega$ is $c_0$-LLC-1, $x$, $y\in B\left(x, \frac{6}{5}d(x,y)\right)$, there exists a curve $\alpha\subset B\left(x, \frac{6}{5}c_0d(x,y)\right)$ joining $x$ and $y$. 
	Let $x_0=x$
	and let $x_1$ be the last point on $\alpha$ along the direction from $x$ to $y$ such that for some $\gamma_{1}\in\Lambda_{xx_1}(\Omega)$, it holds $$\ell(\gamma_{1})\leq d(x, y).$$
	Based on \eqref{H25-0525-3}, we may repeat the above procedure to find $[\kappa_1]$ successive points $x_i\in \alpha$ such that
	for each $i\in\{1,\cdots, [\kappa_1]\}$, $x_i$ is the last point in $\alpha[x_{i-1},y]$ along the direction from $x_{i-1}$ to $y$
	which satisfies for some $\gamma_{i}\in\Lambda_{xx_i}(\Omega)$, it holds
	\be\label{H25-0520-5}
	\ell(\gamma_{i})\leq e^{10\kappa_0\kappa_2i}d(x, y),
	\ee
	and so
	\be\label{H25-0520-6}
	\ell(\gamma_{i})> e^{10\kappa_0\kappa_2(i-1)}d(x, y)\ge e^{10\kappa_0\kappa_2}d(x, y)
	\ee
	for $i\in\{2,\cdots, [\kappa_1]\}$.
	
	For $i\in\{2,\cdots, [\kappa_1]\}$, by \eqref{H25-0520-6}, we may take $x_{0,i}\in \gamma_{i}$ to be the first point along the direction from $x_i$ to $x$ such that
	\be\label{H25-0520-7}
	\ell(\gamma_{i}[x_i,x_{0,i}])=d(x, y).
	\ee
	Then by \eqref{H25-0520-6} and \eqref{H25-0520-7}, it holds
	\be\label{H25-0520-11}
	\ell(\gamma_{i}[x,x_{0,i}])=\ell(\gamma_{i})-\ell(\gamma_{i}[x_i,x_{0,i}])>(e^{10\kappa_0\kappa_2}-1)d(x, y)>\ell(\gamma_{i}[x_i,x_{0,i}]),
	\ee
	which, together with Lemma \ref{Theorem-1.8}, shows that
	\be\label{H25-0520-9}
	d(x, y)=\min\{\ell(\gamma_{i}[x_i,x_{0,i}]),\ell(\gamma_{i}[x,x_{0,i}])\}\leq \kappa_0d_{\Omega}(x_{0,i}).\ee
	
	Let $B_{0,0}=B(x, \frac{6}{5}(c_0+1)d(x, y))$. For each $i\in \{2, \dots, [\kappa_1]\}$, set
	$$B_{1,i}=B\left(x_{0,i}, \frac{1}{5\kappa_0}d(x, y))\right).$$
	Then {\eqref{H25-0520-9} implies that $\overline{B_{1,i}}\subset B\left(x_{0,i}, \frac{1}{2}d_{\Omega}(x_{0,i})\right) \subset \Omega$}.
	Since $\alpha\subset B\left(x, \frac{6}{5}c_0d(x,y)\right)$, for each $u \in \overline{B_{1,i}}$, we have
	$$ d(x, u)\leq d(x,x_i)+d(x_i,x_{0,i})+d(x_{0,i},u)\stackrel{\eqref{H25-0520-7}}{\leq}(\frac{6}{5}c_0+1+\frac{1}{5\kappa_0})d(x,y)<\frac{6}{5}(c_0+1)d(x,y),
	$$
	which implies $\overline{B_{1,i}} \subset B_{0,0}$.
	Then by $[\kappa_1] - 2>2Q^{[\log_212\kappa_0c_0]}$, Lemma \ref{qs-5} (with $R=\frac{6}{5}(c_0+1)d(x, y)$ and $r=\frac{1}{5\kappa_0}d(x, y)$) implies that there exist two integers $s,t$ satisfying $s+1< t \in \{3, \dots, [\kappa_1] \}$ such that
	$$\overline{B_{1,s}}\cap \overline{B_{1,t}}\not=\emptyset,$$
	and thus
	\[
	d(x_{0,t},x_{0,s})\leq \frac{2}{5\kappa_0}d(x,y)\stackrel{\eqref{H25-0520-9}}{\leq} \frac{2}{5} d_{\Omega}(x_{0,t}).
	\]
	Applying Lemma \ref{Lemma-2.1} with $a=\frac{5}{2}$ gives
	$$k_{\Omega}(x_{0,s},x_{0,t})\leq \frac{20}{27}<1,$$
	By \eqref{(2.1)}, \eqref{H25-0525-2} and \eqref{H25-0520-11}, shows that
	\begin{align} \nonumber
		k_\Omega(x,x_{0,s})&\ge\log\left(1+\frac{\ell(\gamma_s[x, x_{0,s}])}{\min\{d_\Omega(x),d_\Omega(x_{0,s})\}}\right)\\
		\nonumber &>\log\left(1+\frac{(e^{10\kappa_0 \kappa_2 }-1)d(x,y)}{\min\{d_\Omega(x),d_\Omega(x_{0,s})\}}\right)\\
		\nonumber &\geq\log\left(1+\frac{e^{10\kappa_0 \kappa_2 }-1}{2}\right)>1.
	\end{align}
	Similarly, we also have
	$$k_\Omega(x,x_{0,t})>1.$$
	Appling {Lemma \ref{Theorem-1.8-0}} with $x=x$, $y=x_{0,t}$ and $z=x_{0,s}$ gives
	$$
	\ell(\gamma_t[x, x_{0,t}]) \le e^{9\kappa_0 \kappa_2} \ell(\gamma_s[x, x_{0,s}]).
	$$
	It follows from the above estimate and \eqref{H25-0520-7} that
	\begin{align} \nonumber
		\ell(\gamma_t) &=\ell(\gamma_t[x, x_{0,t}]) +\ell(\gamma_t[x_t, x_{0,t}])=\ell(\gamma_t[x, x_{0,t}]) +\ell(\gamma_t[x_s, x_{0,s}])\\
		&\le e^{9\kappa_0 \kappa_2}\ell(\gamma_s[x, x_{0,s}]) + \ell(\gamma_s[x_s, x_{0,s}])< e^{9\kappa_0 \kappa_2} \ell(\gamma_s).
	\end{align}
	On the other hand, since $s <t-1$, \eqref{H25-0520-5} and \eqref{H25-0520-6} imply
	$$
	\ell(\gamma_t) > e^{10\kappa_0 \kappa_2 (t-1)} d(x,y) \ge e^{10s\kappa_0 \kappa_2 +10\kappa_0 \kappa_2} d(x,y) \ge e^{10\kappa_0 \kappa_2}\ell(\gamma_s).
	$$
	The above two estimates clearly contradict with each other. Consequently, it holds
	$$\ell(\gamma_{xy})\leq e^{19\kappa_0\kappa_1\kappa_2}d(x,y).$$

\end{proof}

The proof of Theorem \ref{thm:main result uniform} extends to Theorem \ref{thm:main result inner uniform} with only minor differences. 
\begin{proof}[Sketch on the proof of Theorem \ref{thm:main result inner uniform}]
The implications ``(1) $\Rightarrow$ (2)" and ``(2) $\Rightarrow$ (3)" are similar to the previous case and thus is omitted here. We only consider the implication ``(3) $\Rightarrow$ (1)" and outline the differences.  

As in the  of Theorem \ref{thm:main result uniform}, we shall prove that
\begin{itemize}
	\item[i)] For $x,y \in \Omega$, $\gamma_{xy}\in \Lambda_{xy}(\Omega)$, $\alpha_{xy}\in \varGamma_{xy}(\Omega)$, it holds
	\be\label{H25-0520-1-1}
	\ell(\gamma_{xy})\leq e^{19\kappa_0\kappa_1\kappa_2}\ell{(\alpha_{xy})}.
	\ee
	
	\item[ii)] For each $z\in\gamma_{xy}$, it holds \be\label{H25-0520-2-2}\min\left\{\ell(\gamma_{xy}[x,z]),\ell(\gamma_{xy}[y,z])\right\}\leq \kappa_0 d_{\Omega}(z),\ee
	where $\kappa_0=(28ec_0c_1)^2([Q^{2\log_2 70ec_0c_1}]+1)$, $\kappa_1=\kappa_0([Q^{\log_2 15\kappa_0}]+1)$ and $\kappa_2=\kappa_1^{4\kappa_1^2+1}$.
\end{itemize}

	As before, we shall prove \eqref{H25-0520-1-1} via a contradiction argument. 
	Suppose on the contrary that there exists some $\gamma_0\in\Lambda_{xy}(\Omega)$  such that
	$$\ell(\gamma_0)> e^{19\kappa_0\kappa_1\kappa_2}\ell(\alpha_{xy}).$$
	
	Similar to Theorem \ref{thm:main result uniform}, let $x_0=x$. Then we can find  $[\kappa_1]$ successive points $x_i\in \alpha_{xy}$ such that
	for each $i\in\{1,\cdots, [\kappa_1]\}$, $x_i$ is the last point in $\alpha[x_{i-1},y]$ along the direction from $x_{i-1}$ to $y$
	which satisfies for some $\gamma_{i}\in\Lambda_{xx_i}(\Omega)$, it holds
	\be\label{H25-0520-3-2}
	\ell(\gamma_{i})\leq e^{10\kappa_0\kappa_2i}\ell(\alpha_{xy}).
	\ee
	Take $x_{0,i}\in \gamma_{i}$ to be the first point along the direction from $x_i$ to $x$ such that
	\be\label{H25-0520-4-2}
	\ell(\gamma_{i}[x_i,x_{0,i}])=\ell(\alpha_{xy}).
	\ee
	
	Let $B_{0,0}=B(x, 3\ell(\alpha_{xy}))$. For each $i\in \{2, \dots, [\kappa_1]\}$, set
	$$B_{1,i}=B\left(x_{0,i}, \frac{1}{5\kappa_0}\ell(\alpha_{xy})\right).$$ 
	Then $\overline{B_{1,i}}\subset B_{0,0}.$ After that we may argue as in Theorem \ref{thm:main result inner uniform} to derive a contradiction. 
	Therefore, it holds
	$$\ell(\gamma)\leq e^{19\kappa_0\kappa_1\kappa_2}\ell(\alpha_{xy}).$$
\end{proof}


\appendix

\section{Proofs of auxiliary results}

\bpf[Proof of Lemma \ref{lem-3-1}]
We first show that for each $x\in \alpha$, it holds
\be\label{eq:basic on k} 
d_\Omega(x)\geq \max \left\{\frac{2\ell(\alpha[u,x])+d_\Omega(x)}{4c},\frac{1}{2c}d_\Omega(u)\right\}.
\ee 

Indeed, if $\ell(\alpha[u,x])\geq \frac{1}{2}d_{\Omega}(u)$, then $$d_\Omega(u)\leq 2\ell(\alpha[u,x])\leq 2c d_\Omega(x)$$
and so 
\[
\max \left\{\frac{2\ell(\alpha[u,x])+d_\Omega(u)}{4c},\frac{1}{2c}d_\Omega(u)\right\}=\frac{2\ell(\alpha[u,x])+d_\Omega(u)}{4c}\leq d_\Omega(x). 
\] 

If $\ell(\alpha[u,x])< \frac{1}{2}d_{\Omega}(u)$, then $$d_\Omega(x)\geq d_\Omega(u)-\ell(\alpha[u,x])>\frac{1}{2}d_\Omega(u).$$
Since $c\geq 1$, we infer from the above estimate that
$$d_\Omega(x)\geq \frac{1}{2c}d_\Omega(u)=\max \left\{\frac{2\ell(\alpha[u,x])+d_\Omega(u)}{4c},\frac{1}{2c}d_\Omega(u)\right\}.$$

In either case, we proved \eqref{eq:basic on k}.

For each $w\in \alpha$, by \eqref{eq:basic on k}, we have 
\beqq
\begin{aligned}
	k_{\Omega}(u,w)&\leq \ell_k(\alpha[u,w])=\int_{x\in\alpha[u,w]}\frac{ds}{d_{\Omega}(x)}\leq 4c\int_{x\in\alpha[u,w]}\frac{ds}{2\ell(\alpha[u,x])+d_\Omega(u)}\\ &\leq 2c \log\Big(1+\frac{2\ell(\alpha[u,w])}{d_\Omega(u)}\Big)\leq 4c \log\Big(1+\frac{\ell(\alpha[u,w])}{d_\Omega(u)}\Big).
\end{aligned}
\eeqq
\epf

\begin{proof}[Proof of Lemma \ref{Lemma-2.1} ]
	Since $(X,d)$ is a length space,  for each $\varepsilon\in (0,(9a+1)^{-1})$, there exists some curve $\alpha=\alpha_{\varepsilon}$ in $X$ connecting $x_1$ and $x_2$ such that
	\beq\label{20-7-25-1} \ell(\alpha)\leq \big(1+(a-1)\varepsilon\big)d(x_1,x_2).
	\eeq

	We claim that $\alpha\subset \Omega.$ Indeed, if not, then there exists some point $z\in \alpha\cap \partial \Omega$ and thus it follows from \eqref{20-7-25-1} and the assumption $d(x_1,x_2)\leq a^{-1}d_{\Omega}(x_1)$ that
	$$d_{\Omega}(x_1)\leq  d(z,x_1)\leq \ell(\alpha)\leq \big(1+(a-1)\varepsilon\big)d(x_1,x_2)\leq \frac{\big(1+(a-1)\varepsilon\big)}{a}d_{\Omega}(x_1).$$
	It follows that $\varepsilon\geq 1$, which clearly contradicts with our choice of $\varepsilon$.
	
	Let $x\in \alpha$. Since $d(x_1,x_2)\leq a^{-1}d_{\Omega}(x_1)$, \eqref{20-7-25-1} gives
	$$
	d_{\Omega}(x) \geq d_{\Omega}(x_1)-\ell(\alpha)\geq a^{-1}(a-1)(1-\varepsilon)d_{\Omega}(x_1).
	$$
	Then it follows from the above estimate and our assumption $d(x_1,x_2)\leq a^{-1}d_{\Omega}(x_1)$ that
	\beq\label{20-7-25-2}
	\begin{aligned}
		\log \Big(1+\frac{\ell([x_1,x_2])}{d_{\Omega}(x_1)}\Big) &\stackrel{\eqref{(2.1)}}{\leq}  k_{\Omega}(x_1,x_2)\leq \int_{\alpha}\frac{|dx|}{d_{\Omega}(x)}\\
		&\leq \int_{\alpha}\frac{a|dx|}{(a-1)(1-\varepsilon)d_{\Omega}(x_1)}
		\stackrel{\eqref{20-7-25-1}}{\leq} \frac{a\big(1+(a-1)\varepsilon\big)}{(a-1)(1-\varepsilon)}\cdot\frac{d(x_1,x_2)}{d_{\Omega}(x_1)}\\
		&\leq \frac{10a}{9(a-1)}\cdot \frac{d(x_1,x_2)}{d_{\Omega}(x_1)}\leq \frac{10}{9}(a-1)^{-1}.
	\end{aligned}	
	\eeq
\end{proof}


\bpf[Proof of Lemma \ref{Theorem-1.8}]
Select $x_0 \in \gamma$ such that
\beq\label{222}
\min \left\{\diam( \gamma[x_1, x_0]), \diam( \gamma[x_2, x_0])\right\} \geq \frac{1}{3} \diam(\gamma),
\eeq

We first prove that for each $i \in \{1,2\}$ and $x\in  \gamma[x_i, x_0]$,
\be\label{H25-0729-1}
\diam( \gamma[x_i ,x])\leq 10c_0 c_1d_{\Omega}(x).
\ee
Without loss of generality, we only consider the case $x\in \gamma[x_1,x_0]$
and shall prove that 
\beq\label{333}
\diam(\gamma[x_1, x]) \leq 10c_0 c_1d_{\Omega}(x).
\eeq

Let $z_1\in \gamma([x_1,x])\cap \Omega \backslash\overline{B \left(x,\frac{1}{9}\diam( \gamma[x_1, x])\right)}$. It follows from $\diam(\gamma[x_2,x])\geq \frac{1}{3}\diam(\gamma)$ that  there exists $z_2 \in \gamma[x_2,x]$ that belongs to $\Omega \backslash\overline{B \left(x,\frac{1}{9}\diam( \gamma[x_1, x])\right)}$.

Since $\Omega$ is $c_0$-LLC-2, there exists some curve
$$\beta \subset \Omega\backslash \overline{B\left(x,\frac{1}{9c_0}\diam( \gamma[x_1, x])\right)}$$ joining $z_1$ and $z_2$,
and since $(\Omega,\sigma)$ satisfies the $c_1$-ball separation condition, we have
\beqq
B_{\sigma}(x,c_1d_{\Omega}(x)) \cap \beta \neq \emptyset,
\eeqq
from which it follows
\beqq
\diam(\gamma[x_1, x]) \leq 9c_0 c_1d_{\Omega}(x).
\eeqq
This establishes \eqref{333} and thus completes the proof of \eqref{H25-0729-1}.

Next, our aim is to prove that for each $i=1,2$ and $x \in \gamma[x_i, x_0]$, it holds 
\beq\label{444}
\ell (\gamma[x_i, x])\leq (28ec_0c_1)^2([Q^{2\log_2 70ec_0c_1}]+1) d_{\Omega}(x).
\eeq

Without loss of generality, by (\ref{H25-0729-1}), we only consider the case $x\in \gamma[x_1,x_0]$. Take $y_1 \in \gamma[x_1, x]$ such that $$d_{\Omega}(y_1)\geq \frac{1}{2}\sup_{y \in \gamma[x_1, x]}\{d_{\Omega}(y)\}.$$
Then we have
\beq\label{555}
d_{\Omega}(y_1) \leq d_{\Omega}(x)+ d(y_1 ,x)\leq d_{\Omega}(x)+\diam(\gamma[x_1, x]) \stackrel{\eqref{333}}{\leq} (1+10c_0 c_1)d_{\Omega}(x).
\eeq
Note that by \eqref{333}, for each $z\in \gamma[x_1 ,x])$, it holds
\beqq
d(x_1,z) \leq \diam( \gamma[x_1, z]) \leq 10c_0 c_1 d_{\Omega}(z).
\eeqq
Thus we infer from \cite[Lemma 4.4]{GHW} (with $\beta=\gamma[x_1, x]$, $\mu_1=10c_0c_1$ and $x_0=y_1$) that
\beqq
\ell(\gamma[x_1 ,x]) \leq {25}e^3c_0c_1([Q^{2\log_2 70ec_0c_1}]+1)d_{\Omega}(y_1) \stackrel{\eqref{555}}{\leq} ({28}ec_0c_1)^2([Q^{2\log_2 70ec_0c_1}]+1)d_{\Omega}(x).
\eeqq
This establishes \eqref{444} as desired. The proof of Lemma \ref{Theorem-1.8} is thus complete.
\epf


\medskip
\textbf{Acknowledgment.} We would like to thank Prof.~Pekka Koskela,  Dr.~Abhishek Pandey, Yi Xuan, Zhiqiang Yang and Linhong Wang for numerous helpful discussions and comments. 

\end{document}